\documentclass[12pt]{amsart}
\usepackage{amscd,amssymb}
\usepackage[graph,frame,poly,arc]{xy}  
\usepackage[plainpages,backref,urlcolor=blue]{hyperref}
\usepackage{mathtools}

\theoremstyle{plain}
\newtheorem{thm}[subsection]{Theorem}

\newtheorem{prop}[subsection]{Proposition}
\newtheorem{cor}[subsection]{Corollary}

\theoremstyle{definition}
\newtheorem{rk}[subsection]{Remark}

\newtheorem{ex}[subsection]{Example}
\newtheorem{conj}[subsection]{Conjecture}

\numberwithin{equation}{section}
\newcommand{\OO}{{\mathcal O}}

\newcommand{\A}{{\mathcal A}}

\newcommand{\C}{\mathbb{C}}
\newcommand{\PP}{\mathbb{P}}

\begin{document}

\title [Curves with Jacobian syzygies of the same degree]
{Curves with Jacobian syzygies of the same degree}

\author[Alexandru Dimca]{Alexandru Dimca}
\address{Universit\'e C\^ ote d'Azur, CNRS, LJAD, France and Simion Stoilow Institute of Mathematics,
P.O. Box 1-764, RO-014700 Bucharest, Romania}
\email{dimca@unice.fr}

\author[Gabriel Sticlaru]{Gabriel Sticlaru}
\address{Faculty of Mathematics and Informatics,
Ovidius University
Bd. Mamaia 124, 900527 Constanta,
Romania}
\email{gabriel.sticlaru@gmail.com }

\subjclass[2010]{Primary 14H50; Secondary  13D02}

\keywords{plane curve, Milnor algebra, minimal resolution, Tjurina number, line arrangement}

\begin{abstract} 
In this notes we study complex projective plane curves whose graded module of Jacobian syzygies is generated by its minimal degree component. Examples of such curves include the smooth curves as well as the maximal Tjurina curves. However, this class of curves seems to be surprisingly large. In particular, any line arrangement $\A$ of $d=2k+1\geq 5 $ lines having only double and triple points is in this class if the number  of triple points is $k$ and if they are all situated on a line $L \in \A$, see Theorem \ref{thm1X}. 
\end{abstract}
 
\maketitle


\section{Introduction} 

Let $S=\C[x,y,z]$ be the polynomial ring in three variables $x,y,z$ with complex coefficients, and let $C:f=0$ be a reduced curve of degree $d\geq 3$ in the complex projective plane $\PP^2$. 
We denote by $J_f$ the Jacobian ideal of $f$, i.e. the homogeneous ideal in $S$ spanned by the partial derivatives $f_x,f_y,f_z$ of $f$, and  by $M(f)=S/J_f$ the corresponding graded quotient ring, called the Jacobian (or Milnor) algebra of $f$.
Consider the graded $S$-module of Jacobian syzygies of $f$ or, equivalently, the module of derivations killing $f$, namely
\begin{equation}
\label{eqD0}
D_0(f)= \{\theta \in Der(S) \ : \ \theta(f)=0\}= \{\rho=(a,b,c) \in S^3 \ : \ af_x+bf_y+cf_z=0\}.
\end{equation}
This module is also denoted in the literature by 
$AR(f)$ (i.e. all Jacobian relations for $f$) or by ${\rm Syz}(f)$ (i.e. the Jacobian syzygies of $f$). 
We say that $C:f=0$ is an {\it $m$-syzygy curve} if  the module $D_0(f)$ is minimally generated by $m$ homogeneous syzygies, say $\rho_1,\rho_2,...,\rho_m$, of degrees $d_j=\deg \rho_j$ ordered such that $$d_1\leq d_2 \leq ...\leq d_m.$$ 
In this note we assume  that $C$ is not the union of $d$ lines passing through one point, hence $d_1 \geq 1$.
We call these degrees $(d_1, \ldots, d_m)$ the {\it exponents} of the curve $C$ and $\rho_1,...,\rho_m$ a {\it minimal set of generators } for the module  $D_0(f)$. 
The smallest degree $d_1$ is sometimes denoted by ${\rm mdr}(f)$ and is called the minimal degree of a Jacobian relation for $f$. 
It is known that
\begin{equation}
\label{inq1}
2 \leq m \leq d_1+d_2-d+3,
\end{equation}
see \cite[Proposition 2.1]{maxTjurina} and \cite{E2}.

The curve $C$ is {\it free} when $m=2$, since then  $D_0(f)$ is a free module of rank 2, see for instance \cite{KS,Sim2,ST,To}.
A free curve $C$ is also characterized by the condition $d_1+d_2=d-1$, and hence in this case
$r=d_1 <d/2$, see \cite{ST}. 

Recall that the total Tjurina number $\tau(C)$ of a reduced plane curve $C:f=0$ is just the degree of its Jacobian ideal $J_f$, or equivalently the sum of the Tjurina numbers of all the singularities of $C$.
The main result in \cite{duPCTC}, restated below in
Theorem \ref{thmCTC}, gives an upper  bound for the total Tjurina number $\tau(C)$ of $C$ as a function of $d$ and $r$. Moreover, we know that when $r<d/2$ the upper bound is obtained exactly for a free curve, see Corollary \ref{corCTC} below. Moreover, it was shown in \cite{expo} that for any $r<d/2$ there is a free curve with exponents $d_1=r$ and $d_2=d-1-r$. 

The situation of the equality
\begin{equation} 
\label{eqMAX}
\tau(C) = \tau(d,r)'_{max}= (d-1)^2-r(d-r-1)- \binom{2r+2-d}{2}
\end{equation}
for 
$r \geq d/2$ is more subtle.  A reduced curve $C:f=0$
of degree $d$, such that $r={\rm mdr} (f) \geq d/2$ and the equality \eqref{eqMAX} holds for $C$ is called a {\it maximal Tjurina curve of type} $(d,r)$.
In the paper \cite{maxTjurina} we have put forth the following.
\begin{conj}
\label{conjMAX}
For any integer $d\geq 3$ and for any integer $r$ such that $d/2\leq  r \leq d-1 $, there
are maximal Tjurina curves of type $(d, r)$. Moreover, for $d/2\leq  r \leq d -2$, there are maximal
Tjurina line arrangements of type $(d, r)$.
\end{conj}
This conjecture was shown to hold for many pairs $(d,r)$, and in particular for all pairs as above with $d \leq 11$, see \cite{ADS,maxTjurina}, but the general case is still open. In many cases, these constructed maximal Tjurina curves are either line arrangements,
or irreducible curves with interesting geometrical properties, e.g. for the pair $(d,r)=(d,d-1)$ we get as maximal Tjurina curves the maximal nodal curves, namely irreducible nodal curves of degree $d$
having 
$$g=(d-1)(d-2)/2$$
 nodes, see \cite{maxTjurina, Oka}. 
 
 One has the following characterization of maximal Tjurina curves, see \cite[Theorem 3.1]{maxTjurina}.
 \begin{thm}
\label{thmMT}
With the above notation, the curve $C:f=0$ is maximal Tjurina of type $(d,r)$ if and only if
$$r=d_1= d_2= \ldots = d_m \text{ and } m=2r-d+3.$$
\end{thm}
In this sequel note we study curves $C:f=0$ with Jacobian syzygies of the same degree, that is such that
$$r=d_1= d_2= \ldots = d_m,$$
but not necessarily $m=2r-d+3$.
Such a curve is called in the sequel a {\it curve of type} $(d,r,m)$.
Using the inequalities \eqref{inq1}, it follows that for such a curve one has
\begin{equation}
\label{inq2}
2 \leq m \leq m(d,r)_{max}=2r-d+3. 
\end{equation} 
Hence the maximal Tjurina curves are exactly the curves with Jacobian syzygies of the same degree and maximal $m$.
In the sequel we consider the  case of curves of type $(d,r,m)$ satisfying 
$3 \leq m < m(d,r)_{max}$ and set for such a curve
\begin{equation}
\label{Dm}
\Delta m= m(d,r)_{max}-m>0.
\end{equation} 
In this note we show first that a curve $C$ of type $(d,r,m)$ with either 
$m=3$ or $\Delta m =1$ (resp. $\Delta m =2$ or $\Delta m =3$) has the total Tjurina number determined (resp. almost determined) by the pair $(d,r)$ as in the case if maximal Tjurina curves, see Proposition \ref{prop1A}, 
Theorem \ref{thm1B},  Theorem \ref{thm1C} and Theorem \ref{thm1D}. Then we give a number of examples, showing that the curves of type $(d,r,m)$ exist in a lot of cases. In fact we put forth the following.
\begin{conj}
\label{conjT}
For any integers $d,m\geq 3$ and for any integer $r$ such that 
$$\frac{d}{2}\leq  r \leq d-1 \text{ and } m < 2r-d+3,$$
 there
are curves of type $(d,r,m)$.
\end{conj}
The cases when $m=3$ and $d-3 \leq r \leq d-1$ are covered in Section  4.  For each such value of $r$, we give  a countable family of curves of type $(d,r,3)$, see Propositions \ref{prop1}, \ref{prop2} and \ref{prop3}.
 In particular, in Proposition \ref{prop1} we show that a curve $C$ has type $(d,d-1,3)$ if and only if $C$ is a smooth curve of degree $d$.
This section  contains also other examples, with $m=3$ and $r<d-3$.
Section 5 contains many examples of curves of  type $(d,r,m)$ with $m>3$. 

In the last section we recall first some properties of nodal line arrangements, see Proposition \ref{prop1X}. Then we show that  any line arrangement $\A$ of $d=2k+1$ lines  with $k \geq 3$ having only double and triple points is of type $(d,r,m)=(2k+1,2k-2,k)$ if the number of triple points is $k$ and if they are all situated on a line $L \in \A$, see Theorem \ref{thm1X}. This condition on the triple points of $\A$ can be reformulated as follows. Consider any nodal arrangement $\A'$ with an even number of lines $d'=2k\geq 4$. Then by adding a new line $L$ to $\A'$ we get a line arrangement having only double and triple points with
total Tjurina number satisfying
\begin{equation}
\label{TX}
\tau(\A) \leq k(2k+1)+k.
\end{equation} 
If $\tau(\A)$ is maximal, that is if equality holds in \eqref{TX}, then $\A$ is a line arrangement  of type $(d,r,m)=(2k+1,2k-2,k)$. On the other hand, if we start with a nodal arrangement $\A'$ with an odd number of lines, the similar property fails (by a small margin), see Remark \ref{rk1X} for details.

\section{Basic facts and notions related to Jacobian syzygies of curves }

According to Hilbert Syzygy Theorem, the general form of the minimal resolution for the Milnor algebra $M(f)$ of a curve $C:f=0$ that is assumed to be not free, is the following
\begin{equation}
\label{res2A}
0 \to \oplus_{i=1} ^{m-2}S(-e_i) \to \oplus_{j=1} ^mS(1-d-d_j)\to S^3(1-d)  \to S,
\end{equation}
with $e_1\leq ...\leq e_{m-2}$ and $d_1\leq ...\leq d_m$.
It follows from \cite[Lemma 1.1]{HS} that one has
\begin{equation}
\label{res2B}
e_j=d+d_{j+2}-1+\epsilon_j,
\end{equation}
for $j=1,...,m-2$ and some integers $\epsilon_j\geq 1$. Using \cite[Formula (13)]{HS}, it follows that one has
\begin{equation}
\label{res2C}
d_1+d_2=d-1+\sum_{i=1} ^{m-2}\epsilon_j \geq d-1+m-2,
\end{equation}
which gives in particular \eqref{inq1}.
The last morphism in \eqref{res2A} is given by
$$(a,b,c) \mapsto af_x+bf_y+cf_z,$$
and has as kernel precisely the graded $S$-module of Jacobian syzygies $D_0(f)$ shifted by $(1-d)$.

We continue this section by recalling the following result due to  du Plessis and Wall, see \cite[Theorem 3.2]{duPCTC} as well as \cite{E} for an alternative approach.
\begin{thm}
\label{thmCTC}
For positive integers $d$ and $r$, define two new integers by 
$$\tau(d,r)_{min}=(d-1)(d-r-1)  \text{ and } 
\tau(d,r)_{max}= (d-1)^2-r(d-r-1).$$ 
Then, if $C:f=0$ is a reduced curve of degree $d$ in $\PP^2$ and  $r={\rm mdr}(f)$,  one has
$$\tau(d,r)_{min} \leq \tau(C) \leq \tau(d,r)_{max}.$$
Moreover, for $r={\rm mdr}(f) \geq d/2$, the stronger inequality
$\tau(C) \leq \tau'(d,r)_{max}$ holds, where
$$\tau(d,r)'_{max}=\tau(d,r)_{max} - \binom{2r+2-d}{2}.$$

\end{thm}

At the end of the proof of Theorem \ref{thmCTC}, in \cite{duPCTC}, the authors state the following very interesting consequence (of the proof, not of the statement) of Theorem \ref{thmCTC}.
\begin{cor}
\label{corCTC} Let $C:f=0$ be a reduced curve of degree $d$ in $\PP^2$ and  $r={\rm mdr}(f)$. One has
$$ \tau(C) =\tau(d,r)_{max}$$
if and only if $C:f=0$ is a free curve, and then $r <d/2$.
\end{cor}
In the paper \cite{Dmax}, the first author has given an alternative proof of Corollary \ref{corCTC} and has
shown that  a plane curve $C$ is nearly free, which can be defined by the property $m=3$, $d_1+d_2=d$ and $d_2=d_3$ if and only if a similar property holds. Namely, one has the following result, see \cite{Dmax}.

\begin{cor}
\label{corCTC2} Let $C:f=0$ be a reduced curve of degree $d$ in $\PP^2$ and  $r={\rm mdr}(f)$.
One has
$$ \tau(C) =\tau(d,r)_{max}-1$$
if and only if $C:f=0$ is a nearly free curve, and then $r  \leq d/2$.
\end{cor}
The free and nearly free curves have strong relations with the rational cuspidal plane curves, see \cite{DStRIMS,Mosk}. See also \cite{MaVa} for an alternative view on nearly free curves.

\section{Some properties of curves of type $(d,r,m)$}

 \begin{prop}
\label{prop1A}
For a curve $C$ of type $(d,r,m)$, with the above notation, one has
\begin{equation}
\label{I1}
\Delta m=\sum_{i=1} ^{m-2}(\epsilon_j-1).
\end{equation}
In particular, when $m=3$, the type $(d,r,3)$ determines all the degrees in the resolution \eqref{res2A} and 
\begin{equation}
\label{T2}
\tau(C)=3r(d-1-r).
\end{equation} 
 \end{prop}
\proof
The equality \eqref{I1} follows from the equality \eqref{res2C}. When $m=3$, then there is just $\epsilon_1$ to be determined in order to have the degree $e_1$ in the resolution \eqref{res2A}. But the equality \eqref{I1} implies that
$\epsilon_1=\Delta m+1$ and hence 
$$e_1=d+r+ \Delta m =3r+3-m=3r.$$
To get the formula for $\tau(C)$, recall that for a 3-syzygy curve $D$ of degree $d$ with exponents $(d_1,d_2,d_3)$, the corresponding total Tjurina number $\tau(D)$ is given by the formula
\begin{equation}
\label{T1}
\tau(D)=(d-1)(d_1+d_2+d_3)-(d_1d_2+d_2d_3+d_1d_3),
\end{equation} 
see \cite[Proposition 2.1]{minTjurina}.
In particular, if $C$ is a curve of type $(d,r,3)$, we get the equality \eqref{T2}.

\endproof
 \begin{thm}
\label{thm1B}
For a curve $C$  of type $(d,r,m)$ with $m>3$ such that
$\Delta m=1$
one has $\epsilon_{m-2}=2$, $\epsilon_j=1$ for $1 \leq j \leq m-3$ and 
$$\tau(C)=\frac{d(d-1)}{2}-1+r(d-r-2).$$
\end{thm}
\proof
Using the resolution \eqref{res2A}, we see that the Hilbert series of the Milnor algebra $M(f)$ of a curve of type $(d,r,m)$ is given by the fraction
$$A(t)= \frac{1-3t^{d-1}+mt^{d+r-1}-\sum_{j=1}^{m-2}t^{d+r+\epsilon_j-1}}{(1-t)^3}.$$
It is known that $\dim M(f)_j=\tau(C)$ for $j > 3(d-2)$, see for instance
\cite[ Corollary 8]{CD}. It follows that the fraction
$$B(t)=A(t)-\frac{\tau(C)}{1-t}$$
is in fact a polynomial in $t$ and hence has a finite limit when $t \to 1$. As a result, the second derivative of the denominator
\begin{equation}
\label{Den}
D(t)=1-3t^{d-1}+mt^{d+r-1}-\sum_{j=1}^{m-2}t^{d+r+\epsilon_j-1}-\tau(C)(1-t)^2
\end{equation} 
of the fraction $B(t)$, which has a pole of order 3 at $t=1$, has to vanish at $t=1$.
Moreover, the equality $\Delta m=1$ and the equality \eqref{I1} imply that
$\epsilon _{m-2}=2$ and $\epsilon_j=1$ for $1 \leq j \leq m-3$.
It follows that
$$D(t)=1-3t^{d-1}+mt^{d+r-1}-t^{d+r+1} -(m-3)t^{d+r}-\tau(C)(1-t)^2,$$
where $m=2r-d+2$.
The condition $D''(1)=0$ yields the formula for $\tau(C)$, and this completes the proof of Theorem \ref{thm1B}.
\endproof

 \begin{thm}
\label{thm1C}
For a curve $C$ of type $(d,r,m)$ with $m>3$ such that
$\Delta m=2$
one has either
\begin{enumerate}

\item 
$\epsilon_{m-2}=3$,  $\epsilon_j=1$ for $1 \leq j \leq m-3$ and
$$\tau(C)=\frac{d(d-1)}{2}-3+r(d-r-2),$$
or

\item $\epsilon_{m-2}=\epsilon_{m-3}=2$,  $\epsilon_j=1$ for $1 \leq j \leq m-4$ and
$$\tau(C)=\frac{d(d-1)}{2}-2+r(d-r-2).$$

\end{enumerate}

\end{thm}
\proof
The idea of this proof is completely similar to the proof of Theorem \ref{thm1B}.
The condition $\Delta m=2$ and the equality \ref{I1} shows that 
$m=2r-d+1$ and that only the following two cases are possible.

\medskip

\noindent {\bf Case 1.} $\epsilon_{m-2}=3$,  $\epsilon_j=1$ for $1 \leq j \leq m-3$. Then the formula for the denominator $D(t)$ in \eqref{Den} becomes
$$D(t)= 1-3t^{d-1}+mt^{d+r-1}-t^{d+r+2} -(m-3)t^{d+r}-\tau(C)(1-t)^2.$$
The condition $D''(1)=0$ yields now the formula for $\tau(C)$ in Theorem \ref{thm1C} (1). 

\medskip

\noindent {\bf Case 2.} $\epsilon_{m-2}=\epsilon_{m-3}=2$,  $\epsilon_j=1$ for $1 \leq j \leq m-4$. Then the formula for the denominator $D(t)$ in \eqref{Den} becomes
$$D(t)= 1-3t^{d-1}+mt^{d+r-1}-2t^{d+r+1} -(m-4)t^{d+r}-\tau(C)(1-t)^2.$$
The condition $D''(1)=0$ yields in this case the formula for $\tau(C)$ in Theorem \ref{thm1C} (2). 

\endproof

The same proof as for Theorems \ref{thm1B} and \ref{thm1C} yields the following result.

 \begin{thm}
\label{thm1D}
For a curve $C$ of type $(d,r,m)$ with $m>3$ such that
$\Delta m=3$
one has the following cases.
\begin{enumerate}

\item 
$\epsilon_{m-2}=4$,  $\epsilon_j=1$ for $1 \leq j \leq m-3$ and
$$\tau(C)=\frac{d(d-1)}{2}-6+r(d-r-2),$$

\item $\epsilon_{m-2}=3$, $\epsilon_{m-3}=2$,  $\epsilon_j=1$ for $1 \leq j \leq m-4$ and
$$\tau(C)=\frac{d(d-1)}{2}-4+r(d-r-2),$$

\item $\epsilon_{m-2}=\epsilon_{m-3}= \epsilon_{m-4}= 2$,  $\epsilon_j=1$ for $1 \leq j \leq m-5$ and
$$\tau(C)=\frac{d(d-1)}{2}-3+r(d-r-2).$$

\end{enumerate}

\end{thm}

\section{Curves of type $(d,r,3)$}

If we set $m=3$, we get from \eqref{Dm}
$$\Delta m=2r-d>0.$$
On the other hand $r \leq d-1$ for any curve, due to Koszul type syzygies.  
Our first result in this section is the following.
 \begin{prop}
\label{prop1}
The curve $C$ is of type $(d,d-1,3)$ if and only if $C$ is a smooth curve of degree $d$.
\end{prop}
\proof
When $C$ is a smooth curve, then the partial derivatives $f_x,f_y$ and $f_z$ form a regular sequence in $S$ and this implies that 
the curve $C$ is of type $(d,d-1,3)$. To prove the converse, we use the formula \eqref{T2}.
Hence, if $C$ is a curve of type $(d,d-1,3)$, then $\tau(C)=0$ and therefore $C$ is smooth.
\endproof
From now on we look at the existence of curves of type $(d,r,3)$
for 
\begin{equation}
\label{T3}
\frac{d}{2} <r = d-s,
\end{equation} 
with $s \geq 2$. This  forces $d \geq 2s+1$. 
\begin{prop}
\label{prop2}
The curve 
$$C: f= xyz(x^{d-3}+y^{d-3}+z^{d-3})=0$$ is of type $(d,d-2,3)$ for any $d \geq 5$, satisfying $\Delta m=d-4$ and $\tau(C)=3(d-2).$
\end{prop}
\proof
To simplify the notation, we set in this proof $k=d-3$.
The curve $C$ is a nodal curve with 4 irreducible components, namely the 3 lines $x=0$, $y=0$, $z=0$ and the smooth curve $D: x^k+y^k+z^k=0$. Using \cite[Theorem 4.1]{Edin} it follows that
$r=d_1=d_2=d_3=d-2=k+1$. Let $\rho_1,\rho_2$ and $\rho_3$ the first (i.e. lowest degree) three syzygies of $f$.
It is easy to see that one may take
$$\rho_1=(xy^k - xz^k, y^{k+1} + (k+2)yz^k,-(k+2)y^kz - z^{k+1}),$$
$$\rho_2=(0, x^ky + y^{k+1} + (k+1)yz^k,-x^kz - (k+1)y^kz - z^{k+1})$$
and
$$\rho_3=(x^{k+1} + (k+2)xz^k, x^ky - yz^k, -(k+2)x^kz - z^{k+1}).$$
For $j=2,3$, let $g_j$ be the determinant of the $3 \times 3$ matrix with
the first row given by $(x,y,z)$, the second row given by $\rho_1$ and the third row given by $\rho_j$. Then a direct computation shows that
$$
h_2=g_2/f=d(y^k-z^k) \text{ and }
h_3=g_3/f=-d^2z^k.$$
Using now \cite[Corollary 3.4]{minTjurina} we get
$$\tau(C) \geq (d-1)+(d-2)^2-(d-3)^2=3d-6,$$
where the equality holds if and only if $C$ is a 3-syzygy curve.
This is indeed the case, since our curve $C$ has $3k$ nodes on the smooth curve $D$, which are its intersections with the 3 lines $x=0$, $y=0$, $z=0$. There are in addition the vertices of the triangle formed by these 3 lines, hence we get 
$$\tau(C)=3k+3=3(d-3)+3=3d-6.$$
\endproof
\begin{rk}
\label{rkm3}
Consider the curve
$$C:f=xyz(x-z)(x-2z)(x-3z)(y-z)(y-2z)(y-3z)(x+y)(x+y-2z)=0.$$
Then $C$ is a 4-syzygy curve with exponents $(6,6,6,7)$ and degree $d=11$. Moreover $\tau(C)=72$, the same value as for a curve of type
$(11,6,3)$ by the formula \eqref{T2}. With the notation from the end of the proof of Proposition \ref{prop2}, the corresponding polynomial
$h_2$ and $h_3$ have the common factor $(x+y)$, and hence the ideal
$(h_2,h_3)$ does not define a complete intersection as required to apply 
\cite[Corollary 3.4]{minTjurina}. 
\end{rk}
\begin{prop}
\label{prop3}
The curve 
$$C: f=xyz(x^{d-4}y+x^{d-5}z^2+z^{d-5}xy+y^{d-4}z)=0$$ is of type $(d,d-3,3)$ for any $d \geq 8$, satisfying $\Delta m=d-6$ and $\tau(C)=6(d-3)$.
\end{prop}
\proof
To simplify the notation, we set again $k=d-3$. Let $\rho_1,\rho_2$ and $\rho_3$ the first (i.e. lowest degree) three syzygies of $f$.
It is easy to see that one may take
$$\rho_1=(A_1, A_2,A_3),$$
where 
$$A_1= (k-4)xy^{k-1} - 5x^{k-1}z - (2k-4)x^2yz^{k-3},$$
$$A_2= (2k-1)y^k + (2k+1)x^{k-2}yz + (k^2-k-2)xy^2z^{k-3},$$
$$A_3=-(k^2-2)y^{k-1}z + (k-2)x^{k-2}z^2 - (2k-4)xyz^{k-2}.$$
Similarly,
$$\rho_2=(B_1, B_2,B_3)$$
where
$$B_1= (k-4)x^{k-1}y + (3k-2)x^{k-2}z^2 + (k^2-k-4)xyz^{k-2},$$
$$B_2= (2k-1)x^{k-2}y^2 + (2k-5)x^{k-3}yz^2 - (k-5)y^2z^{k-2},$$
$$B_3= -(k^2-2)x^{k-2}yz - (k^2-k-1)x^{k-3}z^3 - (2k-2)yz^{k-1}.$$
And finally one takes
$$\rho_3=(C_1, C_2, C_3),$$
where
$$C_1= 5x^k + (3k-2)xy^{k-2}z + (7-k)x^2z^{k-2},$$
$$C_2= -(2k+1)x^{k-1}y + (2k-5)y^{k-1}z + (k^2-2k-5)xyz^{k-2},$$
$$C_3= -(k-2)x^{k-1}z - (k^2-k-1)y^{k-2}z^2 - (2k-4)xz^{k-1}.$$

For $j=2,3$, let $g_j$ and $h_j$ be as defined in the proof of Proposition \ref{prop2}. Then a direct computation shows that
$$
h_2=g_2/f=-d^2z((2k-3)x^{k-3} + (k-2)^2yz^{k-4})$$
{ and }
$$h_3=g_3/f=  -d^2((2k-3)y^{k-2} + (k-2)xz^{k-3}).$$
Using now \cite[Corollary 3.4]{minTjurina} we get
$$\tau(C) \geq 2(d-1)+(d-3)^2-(d-5)^2=6d-18,$$
where the equality holds if and only if $C$ is a 3-syzygy curve.
This is indeed the case, since our curve $C$ has $4$ irreducible components, namely the 3 lines $C_1: x=0$, $C_2:y=0$ and $C_3:z=0$,
and the irreducible curve
$$C_4: x^{d-4}y+x^{d-5}z^2+z^{d-5}xy+y^{d-4}z=0.$$
The curve $C_4$ has a unique singularity, namely a node $A_1$ situated at the point $p_1=(0:0:1)$. The irreducible components $C_j$
intersects each other at the point $p_1$ and at two other points, namely
$p_2=(1:0:0)$ and $p_3=(0:1:0)$. We have to compute the 3 local Tjurina numbers $\tau(C,p_j)$ for $j=1,2,3$.
We start with the point $p_2$ and use $(y,z)$ as local coordinates at $p_2$. Then the 3-jet of the equation $f$ of $C$ at the point $p_2$ is
$y^2z$, and hence the singularity $(C,p_2)$ is a simple singularity $D_m$, for some $m \geq 4$, see for instance \cite[Proposition (8.17)]{RCS}. This singularity, being weighted homogeneous, satisfies $\mu(C,p_2)=\tau(C,p_2)$.
Recall that for two isolated plane curve singularities $(X,0)$ and $(Y,0)$ with no common component one has 
\begin{equation}
\label{eqMilnor}
\mu(X\cup Y,0)=\mu(X,0)+\mu(Y,0)+2(X,Y)_0-1,
\end{equation}
see \cite[Theorem 6.5.1]{CTC}. Using this result for $(X,0)=(L_3,p_2)$
and $(Y,0)=(C_4,p_2)$, we get $\mu(L_3 \cup C_4,p_2)=1$. Then, using
\eqref{eqMilnor} for $(X,0)=(L_2,p_2)$ and $(Y,0)=(L_3 \cup C_4,p_2)$
we get
$$\tau(C,p_2)=\mu(C,p_2)=6.$$
Consider now the point $p_3$ and use $(x,z)$ as local coordinates at $p_2$. Then the 3-jet of the equation $f$ of $C$ at the point $p_3$ is
$xz^2$, and hence, as above, the singularity $(C,p_3)$ is a simple singularity $D_m$. This singularity, being weighted homogeneous, satisfies $\mu(C,p_3)=\tau(C,p_3)$.
Using \eqref{eqMilnor} exactly as above, we get in this case
$$\tau(C,p_3)=\mu(C,p_3)=2k.$$
Consider finally the point $p_1$ and use $(x,y)$ as local coordinates at $p_1$. Then the 4-jet of the equation $f$ of $C$ at the point $p_1$ is
$x^2y^2$, and  the singularity $(C,p_1)$ is a unimodal singularity $T_{2,p,q}$ for some integers $p \geq q \geq 5$. This singularity, is no longer weighted homogeneous, but it is known to satisfy 
$$\tau(C,p_1)=\mu(C,p_1)-1.$$
Using \eqref{eqMilnor} exactly as above, we get in this case
$$\tau(C,p_1)=\mu(C,p_1)=4k-6.$$
The above computations imply that
$$\tau(C)=6+2k+4k-6=6k=6d-18,$$
which completes the proof of our claim.

\endproof

The claims in the following examples follow by a direct computation, using for instance the software SINGULAR \cite{Sing}.
\begin{ex}
\label{exd-3}
The curves
$$C: f= (x-2y+z)(y-3z+x)(z-5x+y)(7x+y+z)(x+11y+z)(x+y+13z)(x+y+z)=0$$ 
and
$$C': f=y(4y^2-z^2)((x+y)^2-z^2)((x-y)^2-z^2)=0$$
are of type $(d,d-3,3)$ for $d=7$ and $\Delta m=1$. Here $s=3$, hence $d=7$ is the minimal degree where such an example may occur by \eqref{T3}.
\end{ex}

\begin{ex}
\label{exd-4}
The curve
$$C: f= xyz(x^2-2y^2)(x^2-3z^2)(y^2-5z^2)=0$$
is of type $(d,d-4,3)$ for $d=9$ and $\Delta m=1$.
The curve
$$C': f= xyz(x+y+z)(x-2y+z)(y-3z+x)(z-5x+y)(7x+y+z)(x+11y+z)(x+y+13z)=0$$
is of type $(d,d-4,3)$ for $d=10$ and $\Delta m=2$. Here $s=4$, hence the minimal degree for such examples is $9$.
\end{ex}

\begin{ex}
\label{exd-5}
The curve 
$$C:f=xyz(x-z)(y-z)(x-2z)(y-2z)(x-3z)(y-3z)(y+z-2x)(x+z-2y)=0$$
is of type $(d,d-5,3)$ for $d=11$ and $\Delta m=1$.
The curve 
$$C': f=xyz(x^3+y^2z)(y^3+z^2x)(z^3+x^2y)=0$$
is of type $(d,d-5,3)$ for $d=12$ and $\Delta m=2$.
The curve 
$$C'': f=xyz(x+y+z)(2x+y+z)(3y+z+x)(5z+x+y)(7x+y+z)(11y+x+z)(13z+x+y)$$
$$(17x+y+z)(19y+x+z)(23z+x+y) =0$$ 
is of type $(d,d-5,3)$ for $d=13$ and $\Delta m=3$. Here $s=5$, hence the minimal degree for such examples is $11$.

\end{ex}

\begin{ex}
\label{exAA}
The curve 
$$C: f=z(x^k-z^k)(y^k-z^k)((x+y)^k-2z^k)=0$$ 
is of type $(3k+1,2k,3)$ for $3 \leq k \leq 20$ and $\Delta m=k-1$. Here $s=k+1$.
\end{ex}

\section{Curves of type $(d,r,m)$ with  $m \geq 4$}

The claims in the  examples of this section follow by a direct computation, using for instance the software SINGULAR \cite{Sing}. Our Theorems \ref{thm1B} and \ref{thm1C} can be used to check the total Tjurina number given by SINGULAR and get additional information on the corresponding Milnor algebra resolution \eqref{res2A}. Indeed, one has the following.

\begin{enumerate}

\item The curves $C, C', C''$ and $C_6$ in Example \ref{ex50AA}, the curves $C_7$ and $C'$ in Example \ref{ex50A}, the curves $C$ and $C_8$ in Example \ref{ex50} and the curves $C'$ and $C'''$ in Example \ref{ex51} satisfy
$$\Delta m=1.$$
Hence the corresponding Tjurina numbers may be computed using 
Theorem \ref{thm1B}.

\item The curve $C_7$ in  in Example \ref{ex50AA}, the curve $C_8$ in Example \ref{ex50A}, the curve $C_9$ in Example \ref{ex50} and the curves $C_{10}$, $C$ and $C''$ in Example \ref{ex51}  satisfy
$$\Delta m=2.$$
Hence the corresponding Tjurina numbers can be compared with the values given by
Theorem  \ref{thm1C} and in this way we can determine completely the degrees in the resolution \eqref{res2A}.
In fact, in all these cases  when we have $\Delta m=2$, different from $C_{10}$ and $C$,  the case $(1)$ of  \ref{thm1C} occurs. For the curves  $C_{10}$ and $C$, the case $(2)$ in Theorem \ref{thm1C} is realized.

\item The curve $C_8$ in Example \ref{ex50AA}, the curve $C_9$ in Example \ref{ex50A}, the curves $C'$ and $C_{10}$ in Example \ref{ex50} and the curve $C_{11}$ in Example \ref{ex51}  satisfy
$$\Delta m=3.$$
Hence the same remark as at the previous case applies. The curves  $C_8$ in Example \ref{ex50AA},
$C_9$ in Example \ref{ex50A} and $C_{10}$ in Example \ref{ex50} correspond to the case $(1)$ in Theorem \ref{thm1D}, the curve
$C_{11}$ in Example \ref{ex51} corresponds to the case $(2)$ in this result and the curve $C'$  in Example \ref{ex50} corresponds to the case $(3)$ in Theorem \ref{thm1D}.
\end{enumerate}

\begin{ex}
\label{ex50AA}
The curve 
$$C: f=z(x^2+y^2-z^2)(x+2y-3z)(x^2-yz)=0$$
is of type $(d,d-2,4)$ with $d=6$, $\tau(C)=14$. Note that
$\Delta m=1$ in this case.
The curve 
$$C':f=xyz(x+y-z)(x+y+3z)(y+z)(x+5y-z)(7x+y-z)=0$$
is of type $(d,d-3,4)$ with $d=8$, $\tau(C)=32$. Note that for this curve
$\Delta m=1.$
The curve 
$$C'':f=xyz(x+y+z)(x^2+y^2+z^2)(x^2-yz)(y^2-xz)=0$$
is of type $(d,d-4,4)$ with $d=10$, $\tau(C)=56$. Note that for this curve
again $\Delta m=1.$

The family of curves 
$$C_d: f=(x-2y)(y-3z)(z-5x)(2y+3z+5x)(x^{d-4}+y^{d-4}-z^{d-4})=0$$
is of type $(d,d-2,4)$ with $\tau(C_d)=2(2d-5)$ for all $  6 \leq d \leq 20$. Note that $\Delta m=d-5$ in this case.

\end{ex}

\begin{ex}
\label{ex50A}
The family of curves 
$$C_d: f=z(x^{d-5}+y^{d-5}-z^{d-5})(x^2-yz)(y^2-xz)=0$$
is of type $(d,d-2,5)$ with $\tau(C_d)=5(d-3)$ for all $  7 \leq d \leq 20$. Note that $\Delta m=d-6$ in this case.
The curve 
$$C':f= xyz(x+y-z)(2x-y)(x+y+3z)(y+z)(x+5y-z)(7x+y-z)=0$$
is of type $(d,d-3,5)$ with $d=9$, $\tau(C)=32$. Note that for this curve
$\Delta m=1.$

\end{ex}

\begin{ex}
\label{ex50}
The curve 
$$C: f=x(x^3+y^2z)(y^3+z^2x)(z^3+x^2y)=0$$
is of type $(10,7,6)$ with $\tau(C)=51.$ Note that
$\Delta m=1$ in this case.
The curve 
$$C': f=(x^4+y^2z^2)(y^4+x^2z^2)(z^4+x^2y^2)=0$$
is of type $(12,9,6)$ with $\tau(C)=72.$ Note that $\Delta m=3$ for this curve.
The family of curves
$$C_d:f_d=(x^{d-6}+y^{d-6}-z^{d-6})(yz-x^2)(xz-y^2)(xy-z^2)=0$$
has type $(d,d-2,6)$ for $8 \leq d \leq 20$ and $\tau(C_d)=3(2d-7)$. We have $\Delta m=d-7$ for the curve $C_d$.

\end{ex}

\begin{ex}
\label{ex51}

The family of curves
$$C_d:f_d=(x-z)(2x+3y-5z)((x^2-yz)(y^3+xz^2)(x^{d-7}+y^{d-7}-z^{d-7}=0$$
has type $(d,d-2,7)$ for $8 \leq d \leq 20$ and $\tau(C_d)=8d-37$. We have $\Delta m=d-8$ for the curve $C_d$.
The curve 
$$C: f=(x-z)(y-z)(x-y)(x^2-yz)(y^3 + z^2x)(x^2+y^2-z^2)=0$$
is of type $(10,7,5)$ with $\tau(C)=50.$ Note that in this case $\Delta m=2.$

The curve 
$$C': f=xyz(x+z)(x+y-z)(2x-y)(x+y+3z)(y+z)(x+5y-z)(7x+y-z)=0$$
is of type $(9,7,7)$ with $\tau(C')=35.$ Note that in this case $\Delta m=1.$
The curve 
$$C'': f=(x^3-y^3+z^3)(x^3+y^2z)(y^3+z^2x)(z^3+x^2y)=0$$
is of type $(12,10,9)$ with $\tau(C'')=63.$ Note that
$\Delta m=2$ for this curve.

The curve 
$$C''': f=(x-y+z)(x^2-y^2+z^2)(x^3+y^2z)(y^3+z^2x)(z^3+x^2y)=0$$
is of type $(12,10,10)$ with $\tau(C''')=65.$ Note that $\Delta m=1$ for this curve.

\end{ex}

\section{A geometric example of line arrangements of type $(d,r,m)=(2k+1,2k-2,k)$ for any $k\geq 2$}

The nodal line arrangements satisfy the following.

 \begin{prop}
\label{prop1X}
Let $\A:f=0$ be a line arrangement of $d \geq 4$ lines in $\PP^2$. 
Then the following properties are equivalent

\begin{enumerate}

\item The line arrangement $\A$ has 
only double points, in other words $\A$ is a generic, or a nodal line arrangement.

\item The total number of points of multiplicity at least 2 in $\A$ is given by $\binom{d}{2}$.

\item $\tau(\A)=\binom{d}{2}$.

\item ${\rm mdr}(f)=d-2$.

\item The line arrangement $\A$ is  maximal Tjurina of type $(d,d-2)$.

\end{enumerate}

\end{prop}

\proof The equivalence of the first 3 properties follows from the well known formulas
$$n_2 +n_3\binom{3}{2}+ n_4\binom{4}{2}+\ldots + n_d\binom{d}{2}=\binom{d}{2}$$
and
$$n_2+4n_3+9n_4+ \ldots + (d-1)^2n_d=\tau(\A),$$
where $n_j$ denotes the number of points of multiplicity $j$ in the arrangement $\A$.

The fact that $(1)$ implies $(4)$ follows from instance from
\cite[Theorem 4.1]{Edin}. The converse implication follows from
\cite[Theorem 1.2]{Mich}.  Now the fact that $(5)$ implies $(4)$ follows from definitions, while the fact that $(1)$ implies $(5)$ is just 
\cite[Proposition 5.11]{maxTjurina}.
\endproof

If we consider now line arrangements with double and triple points, the situation becomes much more complicated, as shown by Ziegler's examples of two such line arrangements $\A_1:f_1=0$ and $\A_2:f_2=0$ having the same combinatorics, but
${\rm mdr}(f_1) \ne {\rm mdr}(f_2)$, see for details \cite{ZP}.
There is however a class of line arrangements with double and triple points where this strange phenomenon cannot occur.
as the following result shows.
 \begin{thm}
\label{thm1X}
Let $\A':f'=0$ be a nodal  arrangement of $d' \geq 2$ lines in $\PP^2$ and let $L$ be a new line, not belonging to $\A'$, such that exactly $s\geq 1$ double points of $\A'$ are situated on $L$. Consider the line arrangement $\A=\A' \cup L:f=0$ and set $d=d'+1=\deg(f)$. Then the following hold.

\begin{enumerate}

\item The line arrangement $\A$ has 
$$n_2= \binom{d}{2}-3s$$ double points and 
$$n_3=s\leq \frac{d'}{2}$$ 
triple points as multiple points and hence
$$\tau(\A)=  \binom{d}{2}+s     .$$
Moreover one has ${\rm mdr}(f)=d-3$.

\item If $d'=2k\geq 4$ is an even number and $s=k$, the maximal number of nodes of $\A'$ that may be collinear,  then the line arrangement
$\A$ is a curve of type $(d,r,m)=(2k+1,2k-2,k)$ for any $k\geq 2$.
Moreover one has $\Delta m=k-2$ and the numerical invariants of the minimal resolution \eqref{res2A} of the Milnor algebra $M(f)$ are determined by $k$. In particular, one has
$\epsilon_j=2$ for all $j=1,k-2$.
\end{enumerate}

\end{thm}
The point $(2)$ says that for $k=2$ we get a free line arrangement $\A$ with exponents $(2,2)$, for $k=3$ we get a maximal Tjurina (in fact a nearly free) line arrangement with exponents $(3,3,3)$,
while for $k>3$ we get a line arrangement  $\A:f=0$ such that the second order syzygies of the Jacobian ideal $J_f$ are quadratic. We recall that for a maximal Tjurina curve one has 
$\epsilon_j=1$ for all $j=1,m-2$, see \cite[Theorem 3.1]{maxTjurina}.
\proof
The proof of $(1)$ is easy. Indeed, it is clear that $\A$ has $s$ triple points and no points of multiplicity $>3$. 
The upper bound for $s$ follows by applying Bezout Theorem to the curves $\A'$ and $L$. To get the number $n_2$ of double points we use the well known formula
$$n_2+3n_3=\binom {d}{2}.$$
The formula for the total Tjurina number is also clear, since each node contributes by 1, and each triple point contributes by 4 to this number.
To prove the formula for ${\rm mdr}(f)$, note first that
${\rm mdr}(f')=d'-2$, see for instance 
\cite[Theorem 4.1]{Edin}. Hence, it follows that
$${\rm mdr}(f) \geq {\rm mdr}(f')=d'-2=d-3$$
using for instance \cite[Proposition 2.12]{ADS}.
On the other hand, since $\A$ has at least a point $p$ with multiplicity
$m_p=3$, it follows from \cite[Theorem 1.2]{Mich}
that
$${\rm mdr}(f)\leq d-m_p=d-3.$$

The proof of $(2)$ is much more involved.
Let $I_f$ denote the saturation of the ideal $J_f$ with respect to the maximal ideal ${\bf m}=(x,y,z)$ in $S$ and consider the following  local cohomology group, usually called the (graded) Jacobian module of $f$, 
 $$N(f)=I_f/J_f=H^0_{\bf m}(M(f)).$$
We set $n(f)_k=\dim N(f)_k$ for any integer $k$.  If we set $T=3(d-2)$, then the sequence $n(f)_k$ is symmetric with respect to the middle point $T/2$, that is one has
\begin{equation}
\label{E1}
n(f)_a=n(f)_b
\end{equation}
for any integers $a,b$ satisfying $a+b=T$, see \cite{Se}.
Using \cite[Theorem 6.2]{DIS}, which is a consequence of the main result in \cite{STY}, we see that there is an exact sequence  for any integer $j$ given by
\begin{equation}
\label{E0}
 0 \to D_0(f')_{j-1} \to  D_0(f)_j \to H^0(L,\OO_{L}(j+1-r)) \to 
\end{equation} 
$$ \to N(f_1)_{j+d'-2} \to N(f)_{j+d-1}   \to H^1(L,\OO_{L}(j+1-r)), $$
where $r$ is the number of points in the intersection  $\A' \cap L$ and $N(f')$ is the graded Jacobian module of $f'$. In our case $r=k$.
Since $\A'$ is a nodal arrangement, it follows from \cite[Theorem 3.2]{CDI} that one has
\begin{equation}
\label{E2}
n(f')_j=0 \text{  for } j\leq d'-2=2k-2.
\end{equation}
Using now \eqref{E1}, it follows that
\begin{equation}
\label{E3}
n(f')_j=0 \text{  for } j\geq 3(d'-2)-(d'-2)=2d'-4=4k-4.
\end{equation}
We apply  the exact sequence \eqref{E0} for $j=d-3$ and get
$$\dim D_0(f)_{d-3}=  \dim H^0(L,\OO_{L}(d-2-k))=
\dim H^0(L,\OO_{L}(k-1))=k,$$
in view of \eqref{E3}. Hence $m$, the minimal number of generators of $D_0(f)$ is at least $k$, that is $m \geq k$, and moreover
$$d_1=\ldots =d_k=d-3=2k-2.$$
Using \eqref{res2C} we get
\begin{equation}
\label{E4}
4k-4=d_1+d_2= d-1+ \sum_{j=1}^{m-2}\epsilon_j=2k+\sum_{j=1}^{m-2}\epsilon_j.
\end{equation}

We apply now the exact sequence \eqref{E0} for $j=d-2$ and recall that
$\dim D_0(f')_{d-3}=d'-1=2k-1$, see for instance 
\cite[Theorem 4.1]{Edin}.
We obtain, in view of \eqref{E2}, that
$$\dim D_0(f)_{d-2}= \dim D_0(f')_{d-3} + \dim H^0(L,\OO_{L}(d-1-k))=$$
$$=(2k-1)+\dim H^0(L,\OO_{L}(k))=3k.$$
Hence the following two cases are possible.

\medskip

\noindent {\bf Case 1.} $\epsilon_j \geq 2$ for all $j$. Then we have using \eqref{E4}
$$2k-4=\sum_{j=1}^{m-2}m\epsilon_j \geq 2(m-2)$$
and hence $m \leq k$. It follows that in this case $m=k$ and $\A$ is a line arrangement of type $(d,r,m)=(2k+1,2k-2,k)$ as claimed in point $(2)$, with
 $\epsilon_j=2$ for all $j$.

\medskip

\noindent {\bf Case 2.} Assume now that exactly $ u \geq 1$ of the $\epsilon_j $ are equal to 1 and the remaining $m-2-u$ satisfy  $\epsilon_j  \geq 2$. Then there are $u$ second order syzygies of degree 1, and to have $\dim D_0(f)_{d-2}=3k$ we need $u$ new generators for $D_0(f)$ in degree $d-2$.
On one hand we would have then $m \geq k+u$, and on the other hand
using \eqref{E4} we get
$$2k-4 \geq u+2(m-2-u)$$
which yields
$$m \leq k+ \frac{u}{2},$$
a contradiction. It follows the second case cannot occur, and this completes our proof.

\endproof

 \begin{rk}
\label{rk1X} (i) As seen in the proof above, the minimal number of generators for $D_0(f')$ is $2k-1$, while for $D_0(f)$ it is just $k$. Hence by the addition of a line, the number of generators may drop a lot, even in cases when the degree of these generators is preserved.

(ii) With the notation from Theorem \ref{thm1X}, if $d'=2k-1 \geq 5$ is an odd integer, and $s=k-1$ is again the maximal number of double points of $\A'$ that may be collinear, then for the new arrangement $\A$, the
graded module of Jacobian syzygies is NOT generated by its minimal degree component. To see an example, consider the nodal arrangement
$$\A'_k: f'_k=\prod_{j=1,k}(jx-y-j^2z) \prod_{j=1,k-1}(jx+y-j^2z)=0$$
and the associated line arrangement $\A_k:f_k=yf_k'=0$. Then for
$k=3$ the exponents of $\A_k$ are $(3,3,4)$, hence $\A_3$ is a plus-one generated line arrangement. For $k=4$ (resp. $k=5$) the exponents of $\A_k$ are $(5,5,5,6)$ (resp. $(7,7,7,7,8)$).
Therefore the property of Theorem \ref{thm1X} $(2)$ does not seem to hold for $d'$ odd, even if the failure is minimal.
\end{rk}

\end{document}